\documentclass[12pt,reqno]{amsart}
\usepackage{verbatim}
\usepackage{amscd, amsfonts}
\usepackage{epsfig}

\begin{document}

\numberwithin{equation}{section}
\renewcommand{\theequation}{\thesection.\arabic{equation}}
\setcounter{secnumdepth}{2}

% The following establish some abbreviations for combinations and/or
% commonly used letters/terms.
%
% The following are macros for Greek letters & special fonts.
\newcommand{\Ac}{{\mathcal A}}
\newcommand{\Bc}{{\mathcal B}}
\newcommand{\bt}{\tilde{b}}
\newcommand{\bfF}{\bf{F}}
\newcommand{\bfG}{\bf{G}}
\newcommand{\Dc}{{\mathcal D}}
\newcommand{\Ec}{\mathcal{E}}
\newcommand{\Ectd}{\tilde{\Ec}}
\newcommand{\Fc}{{\mathcal F}}
\newcommand{\Fctd}{\widetilde{\Fc}}
\newcommand{\Gma}{{\Gamma}}
\newcommand{\Gc}{\mathcal{G}}
\newcommand{\Gctd}{\tilde{\Gc}}
\newcommand{\Hc}{{\mathcal H}}
\newcommand{\htil}{\tilde{h}}
\newcommand{\Lap}{\Delta}
\newcommand{\Mcl}{{\mathcal M}}
\newcommand{\N}{\mathbb{N}}
\newcommand{\Nb}{\mbox{$\mathbb N$}}
\newcommand{\nab}{\nabla}
\newcommand{\Om}{\Omega}
\newcommand{\Omb}{\overline{\Omega}}
\newcommand{\pal}{\partial}
\newcommand{\Pc}{\mathcal{P}}
\newcommand{\psit}{\tilde{\psi}}
\newcommand{\Rb}{\overline{\R}}
\newcommand{\R}{\mbox{$\mathbb R$}}
\newcommand{\Sc}{{\mathcal S}}
\newcommand{\sg}{\sigma}
\newcommand{\sgt}{\tilde{\sigma}}
\newcommand{\Uc}{{\mathcal U}}
\newcommand{\Uct}{{\widetilde \Uc}}
\newcommand{\vap}{\varphi}
\newcommand{\Vc}{\mathcal{V}}
\newcommand{\Vctd}{\widetilde{\Vc}}
\newcommand{\Wc}{{\mathcal W}}
\newcommand{\Wctd}{\widetilde{\Wc}}

\newcommand{\barr}{\begin{eqnarray}}
\newcommand{\bc}{\begin{center}}
\newcommand{\beq}{\begin{equation}}
\newcommand{\bpf}{\begin{proof} \quad}
\newcommand{\btm}{\begin{thm}}

\newcommand{\earr}{\end{eqnarray}}
\newcommand{\ec}{\end{center}}
\newcommand{\eeq}{\end{equation}}
\newcommand{\epf}{\end{proof}}
\newcommand{\etm}{\end{thm}}

\newcommand{\deq}{:= }
\newcommand{\deqs}{\ :=\ }

\newcommand{\eqs}{\ =\ }
\newcommand{\geqs}{\ \geq \ }
\newcommand{\leqs}{\ \leq \ }
\newcommand{\lts}{\ < \ }
\newcommand{\mns}{\, - \, }
\newcommand{\pls}{\, + \, }
\newcommand{\plms}{\ + \ }

% The following are standard combinations of letters.

\newcommand{\foral}{\qquad \mbox{for all} \quad }
\newcommand{\foreach}{\qquad \mbox{for each} \quad }
\newcommand{\wrt}{ \mbox{with respect to}  }
\newcommand{\xand}{\quad \mbox{and} \quad }
\newcommand{\xfor}{ \quad \mbox{for} \ }
\newcommand{\xiff}{\ \mbox{if and only if} \ }
\newcommand{\xon}{\qquad \mbox{on} \quad}
\newcommand{\xor}{\qquad \mbox{or} \quad}
\newcommand{\xthen}{\quad \mbox{then} \quad}
\newcommand{\xwhen}{\qquad \mbox{when} \quad}
\newcommand{\xwith}{\qquad \mbox{with} \quad}

% The following are special for this paper.
\newcommand{\ang}[1]{\langle#1\rangle}
\newcommand{\bdy}{\pal \Om}
\newcommand{\bdyp}{\pal \Omp}
\newcommand{\bdyo}{\pal \Omo}
\newcommand{\bdyt}{\pal \Omt}
\newcommand{\bdypo}{\bdyo \times \Omt}
\newcommand{\bdypt}{ \Omo  \times \bdyt}

\newcommand{\Cinfty}{C^{\infty}}

\newcommand{\Dej}{ \Dnu e_j    }
\newcommand{\delk}{\delta_k}
\newcommand{\deloj}{\delta_{1j}}
\newcommand{\deltk}{\delta_{2k}}
\newcommand{\Dnu}{D_{\nu}}
\newcommand{\dsg}{\, d \sigma}

\newcommand{\Eco}{\Ec_1}
\newcommand{\Ecotd}{\tilde{\Eco}}
\newcommand{\Ect}{\Ec_2}
\newcommand{\Ecttd}{\tilde{\Ect}}
\newcommand{\ej}{e_j}
\newcommand{\eoj}{e_{1j}}
\newcommand{\eok}{e_{1k}}
\newcommand{\etj}{e_{2j}}
\newcommand{\etk}{e_{2k}}
\newcommand{\etjk}{e_j \otimes f_k}
\newcommand{\Eten}{\Ec_{\otimes}}

\newcommand{\Fco}{\Fc_1}
\newcommand{\Fct}{\Fc_2}
\newcommand{\foj}{f_{1j}}
\newcommand{\fok}{f_{1k}}
\newcommand{\ftj}{f_{2j}}
\newcommand{\ftk}{f_{2k}}

\newcommand{\ghj}{\hat{g}_{j}}
\newcommand{\ghjk}{\hat{g}_{jk}}
\newcommand{\ghten}{g \otimes h}

\newcommand{\Gco}{\Gc_1}
\newcommand{\Gct}{\Gc_2}
\newcommand{\goj}{g_{1j}}
\newcommand{\gtk}{g_{2k}}

\newcommand{\grads}{\nabla s}
\newcommand{\gradu}{\nabla u}
\newcommand{\gradv}{\nabla v}
\newcommand{\gradx}{\nabla_x}
\newcommand{\grady}{\nabla_y}

\newcommand{\hhk}{\hat{h}_{k}}
\newcommand{\hhjk}{\hat{h}_{jk}}
\newcommand{\Iby}{\int_{\bdy} \; }
\newcommand{\Ibyo}{\int_{\bdyo} \; }
\newcommand{\Ibyt}{\int_{\bdyt} \; }
\newcommand{\Ibyp}{\int_{\bdyp} \; }

\newcommand{\IOm}{\int_{\Om} \; }
\newcommand{\IOmo}{\int_{\Omo} \; }
\newcommand{\IOmt}{\int_{\Omt} \; }
\newcommand{\IOmp}{\int_{\Omp} \; }

\newcommand{\Lamo}{\Lambda_1}
\newcommand{\Lamt}{\Lambda_2}

\newcommand{\lamj}{\lambda_j}
\newcommand{\lamk}{\lambda_k}

\newcommand{\Lapu}{\Delta \, u}
\newcommand{\Lapv}{\Delta \, v}
\newcommand{\lmoj}{\lambda_{1j}}
\newcommand{\lmok}{\lambda_{1k}}
\newcommand{\lmtj}{\lambda_{2j}}
\newcommand{\lmtk}{\lambda_{2k}}

\newcommand{\n}[1]{\left\vert#1\right\vert}
\newcommand{\nm}[1]{\left\Vert#1\right\Vert}

\newcommand{\Omi}{\Om_i}
\newcommand{\Omp}{\Om_p}
\newcommand{\Omo}{\Om_1}
\newcommand{\Omob}{\bar{\Omo}}
\newcommand{\Omot}{\Omo \times \Omt}
\newcommand{\Omt}{\Om_2}
\newcommand{\Omtb}{\bar{\Omt}}

\newcommand{\opal}{\ \oplus_{\pal} \ }
\newcommand{\oplusb}{\ \oplus_b \ }

\newcommand{\Rni}{{\R}^{N_i}}
\newcommand{\Rno}{{\R}^{N_1}}
\newcommand{\Rnt}{{\R}^{N_2}}

\newcommand{\RN}{{\R}^N}
\newcommand{\ra}{\rightarrow}

\newcommand{\Sco}{\Sc_1}
\newcommand{\Sct}{\Sc_2}
\newcommand{\Scott}{\Sco \otimes \Sct}
\newcommand{\Scz}{\Sc_0}

\newcommand{\soj}{s_{1j}}
\newcommand{\stk}{s_{2k}}

\newcommand{\sumo}{\sum_{j=1}^{\infty}}
\newcommand{\sumz}{\sum_{j=0}^{\infty}}

\newcommand{\ubar}{\overline{u}}
\newcommand{\ut}{\tilde{u}}
\newcommand{\uhj}{\hat{u}_j}
\newcommand{\uhz}{\hat{u}_0}
\newcommand{\ujk}{u_{jk}}
\newcommand{\utjk}{\tilde{u}_{jk}}
\newcommand{\uvten}{u \otimes v}

\newcommand{\vbar}{\overline{v}}
\newcommand{\vjk}{v_{jk}}
\newcommand{\vjktd}{\tilde{v}_{jk}}
\newcommand{\vhj}{\hat{v}_j}
\newcommand{\vhk}{\hat{v}_k}
\newcommand{\vhz}{\hat{v}_0}
\newcommand{\vt}{\tilde{v}}

\newcommand{\wjk}{w_{jk}}
\newcommand{\wjktd}{\tilde{w}_{jk}}

% The following are various spaces used
\newcommand{\Hoo}{H_1}
\newcommand{\Hto}{H_2}
\newcommand{\HFten}{H_1 \otimes_F \Hto}
\newcommand{\Hten}{H_1 \otimes \Hto}

%  Spaces of Lebesgue and Sobolev functions on Omega
\newcommand{\Harm}{\Hc(\Om)}
\newcommand{\HLap}{H(\Delta, \Om)}
\newcommand{\Hone}{H^1(\Om)}
\newcommand{\Honep}{H^1(\Omp)}

\newcommand{\Hono}{H^1(\Omo)}
\newcommand{\Hont}{H^1(\Omt)}

\newcommand{\Honet}{\Hone \times \Hone}

\newcommand{\Hzone}{H_0^1(\Om)}
\newcommand{\Hzoo}{H_0^1(\Omo)}
\newcommand{\Hzot}{H_0^1(\Omt)}
\newcommand{\Hzonep}{H_{0}^1(\Omp)}

\newcommand{\Lp}{L^p (\Om)}
\newcommand{\Lq}{L^q (\Om)}
\newcommand{\Lr}{L^r (\Om)}
\newcommand{\Lt}{L^2 (\Om)}

\newcommand{\Ltp}{L^2 (\Omp)}
\newcommand{\Lto}{L^2 (\Omo)}
\newcommand{\Ltt}{L^2 (\Omt)}

% spaces of smooth and Harmonic functions

\newcommand{\Cc}{C_c^1 (\Om)}

% boundary spaces
\newcommand{\Hhby}{H^{1/2}(\bdy)}
\newcommand{\Hsby}{H^s(\bdy)}

\newcommand{\Linby}{L^{\infty} (\bdy, \dsg)}
\newcommand{\Lpb}{L^p (\bdy, \dsg)}

\newcommand{\Ltby}{L^2 (\bdy, \dsg)}
\newcommand{\Ltbby}{L^2 (\bdy, b \dsg)}
\newcommand{\Lttby}{L^2 (\bdy, \dsgt)}
% end of macros

%  The following are the theorem like structures used here.
%
\newtheorem{thm}{Theorem}[section]
\newtheorem{cor}[thm]{Corollary}
\newtheorem{cond}{Condition}
\newtheorem{lem}[thm]{Lemma}
\newtheorem{prop}[thm]{Proposition}

\title[Tensor Products  and Laplacian Eigenfunctions]
{Laplacian Eigenproblems on Product Regions and  Tensor Products of Sobolev Spaces. }
\author[Auchmuty  \&   Rivas]{Giles Auchmuty \, \& \,   M. A. Rivas}

\address{Giles Auchmuty, \quad
Department of Mathematics, University of Houston, Houston, TX 77204-3008, USA}
\email{auchmuty@uh.edu}

\address{Mauricio A. Rivas, \quad 
Department of Mathematics, Wake Forest University, 
%PO Box 7388, 127 Manchester Hall,
Winston-Salem, NC 27109, USA}
\email{rivasma@wfu.edu}

\thanks {The authors gratefully acknowledges research support by NSF award DMS 11008754. 
 }

\date{\today}

\begin{abstract} 
Characterizations of eigenvalues and eigenfunctions of the Laplacian on a product domain  $\Omp := \Omo\times\Omt$ are obtained.
When zero Dirichlet, Robin or Neumann conditions are specified on each factor, then the eigenfunctions on $\Omp$ are precisely the products of the eigenfunctions on the sets $\Omo, \Omt$ separately. 
There is a related result when Steklov boundary conditions are specified on $\Omt$.
These results enable the characterization of $\Honep$ and $\Hzonep$ as tensor products and descriptions of some orthogonal  bases of the spaces. 
A different  characterization of the trace space of $\Honep$ is found.
\end{abstract}
\maketitle

%%%------------------------------------------------------------------------------------------------------------------------------------------------------------------------------
%%%			SECTION:	   1.	Introduction
%%%-------------------------------------------------------------------------------------------------------------------------------------------------------------------------------

\section{Introduction} \label{s1}

This paper treats some questions related to $H^1$-Sobolev spaces on  {\it product regions} $\Omp := \Omo\times\Omt \subset
\R^{N_1} \times \R^{N_2}$.
A particular interest is the representation of such spaces as tensor products of Sobolev spaces on $\Omo, \Omt$.
Such tensor products provide a formalism for the separation of variables methodology used in approximating solutions of 
partial differential equations. 
In particular we treat questions concerning when products of Laplacian eigenfunctions are orthogonal  bases of these 
Hilbert-Sobolev spaces.

The construction of orthogonal bases of various $H^1$-Sobolev spaces in terms of eigenfunctions of the Laplacian on a region 
is standard and justifies many of the approximation methods used in science and engineering.
In particular it enables a reduction of dimension that is often very important for effective approximation schemes.
A natural question is whether {\bf all} the Laplacian eigenfunctions on a product region are products of Laplacian 
eigenfunctions on the individual factors $\Omo, \Omt$.
This is known to hold when zero-Dirichlet boundary conditions are imposed on both factors $\Omo, \Omt$ as discussed in section \ref{s3}.

The analysis here treats weak ($H^1$-)solutions of the eigenproblem.
This avoids many issues about regularity; some of which arise from the fact that  product regions are not $C^1$ at their ``corners"
and allows us to prove results that hold for many of the situations studied in numerical computations. 
The  case of zero-Neumann boundary conditions is studied in Section \ref{s5}.
Various combinations of eigenproblems with zero Dirichlet, Robin and Neumann  boundary conditions on the factors are described 
and analyzed  in Sections \ref{s6} - \ref{s8}.

The situation where Steklov boundary conditions are prescribed on $\Omt$ and zero Dirichlet, Neumann, or Robin boundary conditions hold 
on $\Omo$ is studied in Section \ref{s9}. 
These lead to a two-parameter eigenproblem on $\Omp$. 
Finally in Section \ref{s11}, a tensor product characterization of the trace space of functions on the boundary of $\Omp$ is described.
This provides a different orthogonal basis for  the trace space of $H^1$ to that described by the first author in  \cite{Au4}.
In this paper our results  are confined to the case where $\Om $ is the product or two regions. 
When the region $\Om$ is the product of any finite number of regions, we expect that similar results will hold but the algebra and the notation becomes much more complicated. 

%%%------------------------------------------------------------------------------------------------------------------------------------------------------------------------------
%%%			SECTION:		Definitions and Notation
%%%-----------------------------------------------------------------------------------------------------------------------------------------------------------------------------
\vspace{1em}

% ####### Section 2 ######
% This is line 410 approx
\section{Definitions and Notation.}\label{s2}

A region $\Om$ is a non-empty, connected, open subset of $\RN$. Its closure  is denoted  $\Omb$ and its boundary 
$\bdy \deq \Omb\setminus \Om$. 
A standard assumption about the region is the following. 

\noindent{\bf (B1):}  {\it $\Om$ is a bounded region in  $\RN$ and its boundary  $\bdy$  is the union of a finite number of 
disjoint closed Lipschitz  surfaces; each surface having finite surface area.}

When this holds there is an outward unit normal $\nu$ defined at $\sigma \; a.e.$  point of $\bdy$. 
The  definitions and terminology of Evans and Gariepy \cite{EG}, will be followed except  that $\sigma, \dsg$,  respectively, 
will represent Hausdorff   $(N-1)$-dimensional measure and integration with respect to this measure. 
All functions in this paper will take values in  $\Rb \deq [-\infty,\infty]$ and derivatives $D_j u$ are  taken in the weak sense.

 The real Lebesgue spaces $\Lp$ and $\Lpb, \ 1 \le p \le \infty$ are  defined in the standard manner and have the usual $p$-norms
  denoted by ${\nm{u}_p}$ and  $ {\nm{u}}_{p,\bdy}$.
 
 Let $\Hone$ and $\Hzone$ be the usual real Sobolev space of  functions on $\Om$.  
 $\Hone$ is a  real Hilbert space under  the standard $H^1$-inner product
\beq\label{e2.3}
{[u , v]}_{1, 2}  \deq \IOm \ [u(x) \, v(x) \pls  \gradu(x)  \cdot  \gradv(x)] \  dx.
\eeq
Here  $\gradu := (D_1u, \ldots, D_Nu)$ is the gradient of the function $u$ and the associated norm  may sometimes be 
 denoted   $ {\nm{u}}_{1, 2, \Om}$ when there could be questions about the region. 

This paper will prove various results about Laplacian eigenproblems on {\it product regions}  $\Omp := \Omo\times\Omt$ with 
$\Omi \subset \Rni$.  
Here points in $\Omo$, $\Omt$ are denoted by $x=(x_1, x_2, \ldots, x_{N_1})$ and $y=(y_1, y_2, \ldots, y_{N_2})$, respectively, 
so that points in the product region $\Omp$ are denoted by  $(x, y)$.  
It is well known that the eigenfunctions of the Laplacian on a region $\Om$ may be found that are  bases of various Lebesgue and 
Sobolev spaces of the form $\Lt, \Hzone$ and $\Hone$. 
When $\Om$ is a product region $\Omp$, then under quite general conditions, it is shown that  various classes of  eigenfunctions of the  Laplacian on $\Omp$ are products of the eigenfunctions of the Laplacian on the factors $\Omo, \Omt$ respectively.
A consequence is  that the Hilbert-Sobolev spaces $\Hzonep, \Honep$ are tensor products of the relevant spaces on $\Omo, \Omt$  respectively.

In this paper we shall use various standard results from the calculus of variations and convex analysis.  
When terms are not otherwise defined they should be taken as in Attouch, Buttazzo and Michaille \cite{ABM}.
Background material on such methods may be found in Blanchard and Br\"{u}ning \cite{BB1} or Zeidler \cite{Z2}.
Some of the notation and many of the methods in this paper are derived from the results developed in the recent paper  of Auchmuty \cite{Au4}.

%%%-----------------------------------------------------------------------------------------------------------------------------------------------------------------------------
%%%		SECTION:		Hilbert Tensor Product spaces.
%%%----------------------------------------------------------------------------------------------------------------------------------------------------------------------------
\vspace{1em}                      
\section{Hilbert Tensor Product spaces.}\label{TP}

          Our interest in studying these Laplacian eigenfunction problems is partially motivated by questions regarding possible
 orthonormal bases of various Hilbert-Sobolev spaces on product regions $\Omp := \Omo\times\Omt$.  
 In particular, we will prove results that generalize the well-known result that $\Ltp = \Lto \otimes \Ltt$  with $\otimes$ 
 denoting the tensor product of the spaces for Lebesgue spaces on $\Omp$.
  
      Here some constructions and  results about Hilbert tensor products will be stated for use in the following sections.
  The description is specially tailored for the current situation; in particular the use of dual spaces is avoided so that distributions are not needed.
  Let $H_1 \subset \Lto,  H_2 \subset \Ltt$ be two real separable Hilbert spaces with inner products $\ang{\cdot,\cdot}_{H_1}$, and $\ang{\cdot, \cdot}_{H_2}$ respectively.
   Let $\Ec_1 := \{e_j : j \in J_1\}, \Ec_2 := \{f_k : k \in J_2\}$ be orthonormal bases of $H_1, H_2$ respectively. 
   Without loss of generality  all  functions are assumed to be Borel measurable and  the {\it tensor product} of $e_j$ and $f_k$ is defined to be the function $e_j\otimes f_k$ given by
   \beq \label{TP.1}
   (\etjk)(x,y) \deqs e_j(x) \, f_k(y) \xfor (x,y) \in \Omp. \eeq    
  Then  $\etjk$ is Borel measurable on $\Omo \times \Omt$ with $\|\etjk\|_{2, \Omp} = \| e_j \|_{2, \Omo} \, \| f_k \|_{2, \Omt}$ from Fubini's theorem.
Functions of this ``separated variables"-form will be called {\it dyads}. 

Let $\Eten := \{ \etjk : j \in J_1, k \in J_2 \}$ and define $\HFten $ to be the vector space of all finite linear combinations of functions in $\Eten$. 
Then $g \in \HFten$ implies that 
\beq \label{TP.3}
g(x,y) \eqs \sum_{j \in J_1, k \in J_2} \ \ghjk \, e_j(x) \, f_k(y) . \eeq
with only finitely many $\ghjk$ non-zero.
Define an inner product on $ \HFten$ by 
\beq \label{TP.5}
\ang{g, h}_{\otimes} \deqs \sum_{j,k} \ \ghjk \, \hhjk . \eeq
Then $\Eten$ will be an orthonormal set in $\HFten$ and, from \eqref{TP.3},
\beq \label{TP.7}
\ghjk \eqs \ang{g, \etjk}_{\otimes}. \eeq

The {\it Hilbert tensor product} $\Hten$ is the completion of $\HFten \ \wrt$ the inner product $\ang{\cdot, \cdot}_{\otimes}$.
It is straightforward to verify that a function $g$ is in $\Hten$ if and only if it has a representation of the form \eqref{TP.3} with
$ \sum_{jk} \ghjk^2 \lts \infty $. 
Thus $\Eten$ is an orthonormal basis of $\Hten$ and the Parseval equality implies that \eqref{TP.5} holds for all $g,h \in \Hten$
with  Fourier coefficients defined by \eqref{TP.7}.
Moreover, when $u \in H_1, \, v \in H_2$, then $\uvten \in \Hten$ with Fourier coefficients $\uhj \vhk  = \ang{u, e_j}_{H_1}\ang{v, f_k}_{H_2}$ and 
\beq \label{TP.8}
\| \uvten \|_{\otimes} \eqs \|u \|_{H_1} \, \|v \|_{H_2}  \eeq

Since two Hilbert spaces with orthonormal bases that are in 1-1 correspondence are linearly isomorphic the following holds.
 
\btm \label{TPt1}
Let $H_1, H_2$ be real separable Hilbert spaces as above with orthonormal bases $\Ec_1, \Ec_2$ and H is a real Hilbert space 
with basis $\Eten$.  Then H is linearly isomorphic to $\Hten$. \etm

A simple result about tensor products that will be used later is the following. 
Its proof is standard analysis.

%%%% Second version of Theorem TPt2
\btm \label{TPt2}
Suppose that $H_1 = V_1 \oplus W_1$ and $H_2 = V_2 \oplus W_2$ are orthogonal decompositions of $H_1, H_2$, respectively.
Then
\beq \label{TP.9}
H_1\otimes H_2 \eqs (V_1\otimes V_2)\oplus (V_1 \otimes W_2) \oplus (W_1\otimes V_2) \oplus (W_1\otimes W_2).
\eeq
Moreover, the tensor product of corresponding orthonormal bases of $V_1, V_2$ is an orthonormal basis for $V_1\otimes V_2$, and a similar result holds for the remaining three terms in the orthogonal decomposition of $H_1\otimes H_2$.
\etm

%%%-----------------------------------------------------------------------------------------------------------------------------------------------------------------------------
%%%		SECTION:		Dirichlet Laplacian Eigenproblems on $\Om$ and $\Omp$.
%%%----------------------------------------------------------------------------------------------------------------------------------------------------------------------------
\vspace{1em}                        
\section{Dirichlet Laplacian Eigenproblems on $\Om$ and $\Omp$.}\label{s3}

The {\it Dirichlet Laplacian eigenproblem} on a region $\Om \subset \RN$ is the problem of finding nontrivial solutions
$(\lambda,u) \in \R \times \Hzone$ of the equation
\beq \label{e3.1}
\IOm \gradu \cdot \nabla h \, \, dx  \eqs \lambda \IOm  \,  u \, h \, dx   \foral h\in \Hzone.
\eeq
This is the weak form of the eigenvalue problem of finding nontrivial solutions of the system
\beq
-\Delta u(x)  \eqs  \lambda \,   u(x) \xon \Om  \xwith  u(x) \eqs  0  \xon \bdy.
\eeq
Here $\Delta$ is the Laplacian on the region $\Om$.

The classical result about these problems may be stated as follows 
(see, for instance, Theorem 8.3.2 and Proposition 8.3.1 of Attouch et al. \cite{ABM}).  
Problem \eqref{e3.1} has a sequence $\Lambda:=\{\lambda_j : j\in \N\}$ of strictly positive eigenvalues, repeated according to multiplicity and $\lambda_j \rightarrow \infty$, and an associated sequence of eigenfunctions $\Ec : = \{e_j : j\in \N\}$.  
The sequence $\Ec$ forms an orthonormal basis of $\Lt$, and the sequence $\Ectd : = \{ \lambda_j^{-1/2} e_j  : j\in \N\}$ forms an orthonormal basis for $\Hzone$ with respect to the inner product $[u, v]_{\nabla} : = \IOm \gradu\cdot \gradv \, dx$, so that
\beq \label{e3.3}
\IOm \nabla e_j \cdot \nabla e_j \, \, dx  \eqs \lambda_j .
\eeq

When $\Omp := \Omot$ is a product region, then the individual  Dirichlet  Laplacian eigenproblems on $\Omo, \Omt$ are 
those of finding nontrivial solutions of equation \eqref{e3.1}
with $\Omo, \Omt$ in place of $\Om$.

When $\Om = \Omo$, the sequence of eigenvalues will be denoted $\Lambda_1 : = \{\lmoj : j\in \N\}$ 
and the eigenfunctions $\Eco = \{ \eoj: j \in \N \}$;
$\Eco$  is an orthonormal basis of $\Lto$.   %%% \Ec changed to \Eco
When $\Om = \Omt$,  the sequence of eigenvalues will be denoted $\Lambda_2 : = \{\lmtk : k\in \N\}$ and the eigenfunctions $\Ect := \{ \etk: k \in \N\}$         %%% missing slash \ in etk.
with $\Ect$ an orthonormal basis of $\Ltt$.
 
 Also, the sequences $\Ecotd := \{\lmoj^{-1/2}\eoj : j \in \N\}$ and $\Ecttd :=\{ \lmtk^{-1/2}\etk : k \in \N\}$ will be orthonormal bases of 
 $\Hzoo, \Hzot$ with respect to the inner products
 \beq 
[u, v]_{\nabla, \Omo} : = \IOmo \gradx u\cdot\gradx v \, \, dx \xand [u, v]_{\nabla, \Omt}\deqs \IOmt \grady u\cdot \grady v \, \, dy
 \eeq
where $\gradx, \grady$ are the corresponding gradient operators on the regions $\Omo, \Omt$.

The following result about  product Dirichlet eigenfunctions is a weak version of Theorem 2 in Chapter 11 of Strauss \cite{St}.
It is worth noting that  there are regularity issues for eigenfunctions of the Laplacian on product domains as
the boundary of  $\Omp$ is not a union of $C^1$-manifolds but instead  has ``corners".
So most   results here are confined  to   weak  solutions of the eigenproblems.

\btm \label{T3.1}       
Assume $\Omo, \Omt$ satisfy $(B1)$, and let $\Eco, \Ect$ be the sequences of Dirichlet Laplacian eigenfunctions on $\Omo, \Omt$
  described above. 
Then $u_{jk}  := \eoj \otimes  \etk$ is a  Dirichlet-Laplacian eigenfunction for the problem \eqref{e3.1} 
on the product region $\Omp = \Omot$ corresponding to the eigenvalue $\lmoj + \lmtk$.
\etm 
\bpf
Let  $h\in C_c^1(\Omp)$.     
Then $h(x, \cdot) \in C_c^1(\Omt)$ and $h(\cdot, y) \in C_c^1(\Omo)$ for any  $(x, y)\in \Omp$.  
Substituting $u= u_{jk}$ in the left side of \eqref{e3.1} yields  
\[   \IOmt \IOmo [ \etk(y) ( \gradx \eoj(x) \cdot \gradx h(x,y))  \pls \eoj(x)(\grady \etk(y)\cdot \grady h(x,y)) ] \, dx\, dy  \]
\[  = \IOmt \etk(y)\IOmo  \gradx \eoj(x) \cdot \gradx h(x,y) \, \, dx\, dy   \pls \IOmo \eoj(x) \IOmt \grady \etk(y)\cdot \grady h(x,y) \, dx\, dy  \] 
 from Fubini's theorem.  
This equals
\[   \lmoj \IOmt \etk(y)\IOmo \eoj(x)\, h(x, y) \, \, dx\, dy + \lmtk\IOmo \eoj(x) \IOmt \etk(y) \, h(x, y) \, \, dy\, dx  \]
upon using the eigenequations on $\Omo, \Omt$.
This is the right side of \eqref{e3.1} with $\lmoj + \lmtk$ and $u_{jk}$ in place of $\lambda$ and $u$, respectively.  
A density argument then yields that $\lmoj + \lmtk$ is an eigenvalue of this problem and $u_{jk}$ is an associated eigenfunction.
\epf

Let $\Uc  := \{ \ujk : j, k \in \N \}$.  
One observes that this set is orthonormal in $\Ltp$. 
For each $j,k$ define $\utjk : = (\lmoj + \lmtk)^{-1/2}u_{jk}$ and let $\Uct$ be the family of all such functions; this set is orthonormal in $\Hzonep$.
\btm \label{T3.2}
Assume $\Omo, \Omt$ satisfy $(B1)$, and let $\Uc,  \Uct$ be the sequences of eigenfunctions on $\Omp$ as described above.
Then $\Uc$ is an orthonormal basis of $\Ltp$ and $\Uct$ is an orthonormal basis of $\Hzonep$.
\etm
\bpf
From the observations in the prepceding paragraph these sets are orthonormal. 
Here we show that the sets are maximal.
Suppose there is a nonzero $h\in \Ltp$ such that
\[   \IOmp u_{jk}(x,y) h(x,y) \, \, dx \, dy \, = 0 \foral j, k\in \N.  \]
Fubini's theorem then implies
\[  \IOmo \eoj(x) \IOmt \etk(y) \, h(x, y) \, \, dy \, dx = 0 \foral j\in \N,  \]
 
Since $\Eco$ is an orthonomal basis of $\Lto$, this implies  $\IOmt \etk(y) \, h(x, y) dy =0$ for $a.e.$ $x\in \Omo$. 
Since  this holds for all $k\in \N$, $h(x, y) =0$ for $a.e.$ $(x, y)\in \Omp$ as $\Ect$ is an orthonormal basis of $\Ltt$.  
This contradicts the assumption that $h$ is non-zero so  $\Uc$ is a maximal set in $\Ltp$.  

Suppose that $\lambda$ is an eigenvalue of this problem but  $\lambda$ not equal to any of the $\lmoj + \lmtk$.
The associated eigenfunction $e$ satisfies \eqref{e3.1} for then is orthogonal to each $u_{jk}$ both in $\Hzonep$ and $\Ltp$.  
The first part of the proof shows this is not possible.  
If, on the other hand, $e$ satisfies \eqref{e3.1} for $\lambda$ equal to a fixed $\lmoj +\lmtk$, then the function
\[
\psi (x, y)  \deqs  e(x,y) - \sum_{\lmoj + \lmtk = \lambda} c_{jk}\, u_{jk} \qquad \text{with } c_{jk} \deqs \frac{\IOmp u \, u_{jk}\,\, dxdy}{\IOmp u_{jk}^2 \,\,  dxdy}
\]
would be orthogonal to each $u_{jk}$ both in $\Hzonep$ and $\Ltp$.  
The first part of the proof, in this case, shows $\psi$ is zero $a.e.$ on $\Omp$, 
or that $e$ belongs to the subspace associated to $\lambda = \lmoj + \lmtk$.  
Hence, the only Dirichlet Laplacian eigenfunctions on $\Omp$ are the dyads $u_{jk} = \eoj\otimes \etk$.

Finally, if $h$ is orthogonal to each $\utjk$ in $\Hzonep$, then from  equation \eqref{e3.1} it must be $L^2$-orthogonal to each $\utjk$
so from the first part of this proof we have a contradiction.
Thus $\Uct$ is maximal in $\Hzonep$.
\epf

\begin{cor} \label{cor3.3}
Assume $\Omo, \Omt$ satisfy $(B2)$, and let $\Uc$ be as described above. Then the only eigenvalues of problem \eqref{e3.1} on the region $\Omp = \Omo\times\Omt$ are $\lmoj +\lmtk$ with $j, k\in \N$, and  $\Hzonep = \Hzoo \otimes \Hzot$.
\end{cor}
\bpf 
This follows from Theorem \ref{TPt1} and the above result.
\epf

\vspace{1em}
%%%-----------------------------------------------------------------------------------------------------------------------------------------------------------------------------
%%%		SECTION:	  	Equivalent Inner Products and Boundary Regularity.
%%%----------------------------------------------------------------------------------------------------------------------------------------------------------------------------

\section{Equivalent Inner Products and Boundary Regularity}\label{s4}

 To extend  the preceding results to Robin, Neumann and Steklov eigenproblems some further conditions must be imposed
 on the regions $\Omo, \Omt$.
 These conditions are essentially conditions on the  regularity of the boundaries of the region that yield  compact imbedding results for $\Hone$.
 For simplicity the following definitions are given in terms of a generic region $\Om \subset \RN$  satisfying (B1).

 The region $\Om$ is said to satisfy {\it Rellich's theorem} provided the imbedding of $\Hone$ into $\Lp$ is compact for 
 $1\leq p < p_s$ where $p_s : = 2N/(N - 2)$ when $N \geq 3$, or $p_s = \infty$ when $N=2$.
There are a number of different criteria on $\Om$ and $\bdy$ that imply this result.
When (B1) holds it is Theorem 1 in Section 4.6 of \cite{EG}; see also Amick \cite{Amick}.

When (B1) holds and $u\in W^{1, 1}(\Om)$ then the trace of $u$ on $\bdy$ is well-defined and is a Lebesgue integrable 
function with respect to $\sigma$;  see \cite{EG}, Section 4.2 for details.
The trace map $\gamma$ is the linear extension of the map restricting Lipschitz continuous functions on $\Omb$ to $\bdy$.
In surface integrals, we will often use $u$ in place of $\gamma (u)$ when considering the trace of a function on $\bdy$.
The region $\Om$ is said to satisfy the $L^2$-{\it compact trace theorem} provided the trace mapping $\gamma: \Hone \ra \Ltby$ is compact.  
The regions for our problems will generally be required to satisfy 

%\vspace{1ex}
\noindent   {\bf (B2):} \quad  {\it $\Om$ and $\bdy$ satisfy (B1), the Rellich theorem and the $L^2$-compact trace theorem.}

Many of the following results are based on choosing appropriate inner products on $\Hone$. 
These will include the {\it $b$-inner product} on $\Hone$ given by
\beq \label{e4.2}
[u,v]_b \deqs  \IOm \gradu \cdot \gradv \, dx \pls \Iby b\, u\, v\, \dsg, \
\eeq
The corresponding norm will be denoted $\|u\|_b$.
When there are different regions involved in the same problem these inner products and norms are written as $[u,v]_{b, \Om}$ and $\|u\|_{b, \Om}$, accordingly. 

The boundary integral in \eqref{e4.2} has a weight function $b: \bdy \to (0,\infty)$ that satisfies 

\noindent   
{\bf (B3):} \quad  {\it $b$ is Borel measurable,  $b \in \Linby$ with $b \geq  b_0 > 0 \ \sg$-a.e. on $\bdy$.}

The associated boundary quadratic form is denoted   $\Bc(u) := \Iby b\, u^2 \dsg$, and it is shown in Theorem 3.1 
of \cite{Au3} that $\Bc$ is weakly continuous on $\Hone$.

 The following result about these bilinear forms is crucial for many   later results.
 \btm \label{T4.1}
 Assume $(B1)$-$(B3)$ hold and the bilinear form $[\cdot, \cdot]_b$ is given by \eqref{e4.2}.  Then this $b$-inner product is equivalent to the standard inner product on $\Hone$ given by \eqref{e2.3}.   \etm
 \bpf
The proof of Theorem 7.2 of \cite{Au4} shows \eqref{e4.2} defines an inner product equivalent to the standard $H^1$-inner
product.
 \epf

\vspace{1em}  
%%%-----------------------------------------------------------------------------------------------------------------------------------------------------------------------------
%%%		SECTION: 		Neumann Laplacian Eigenproblems on $\Om$ and $\Omp$.
%%%----------------------------------------------------------------------------------------------------------------------------------------------------------------------------

\section{ Neumann Laplacian Eigenproblems on $\Om$ and $\Omp$.} \label{s5}

The {\it  Neumann Laplacian eigenproblem} on a region $\Om \subset \RN$ is the problem of finding nontrivial solutions
$(\lambda,u) \in \R \times \Hone$ of the system
\beq \label{e5.1}
\IOm \gradu \cdot \nabla h \, \, dx   \eqs \lambda \IOm     u \, h \, dx   \foral h\in \Hone.
\eeq
This is the weak form of the eigenvalue problem of finding nontrivial solutions of the system
\beq
-\Delta u(x)  \eqs  \lambda \,    u(x) \xon \Om  \xwith  \Dnu u  \eqs  0  \xon \bdy
\eeq
with $\Dnu u := \gradu \cdot \nu$.

A first result about these  Neumann Laplacian eigenproblems on bounded regions  is well-known subject to various 
different assumptions about the region $\Om$. 
Here a statement that can be proved using  the analysis of Section 4 of \cite{Au4} will be used. 
Take $V = \Hone, \ H = \Lt$ in that paper and let
\beq \label{e5.3}
a(u,v) \deqs  [u,v]_{1, 2} \xand m(u,v) \deqs \ang{u,v}_2 \eeq
where $[\cdot, \cdot]_{1, 2}$ and $\ang{\cdot, \cdot}_2$ are the standard $H^1$- and $L^2$-inner products, respectively.

Obviously $\lambda = 0$ is an eigenvalue of \eqref{e5.1} with the associated eigenfunction $f(x) \equiv 1$.
To obtain the sequence of all Neumann Laplacian eigenvalues, 
consider the sequence of constrained variational problems $\Pc_k$ as described by  (4.1) in \cite{Au4}.
This generates a sequence $\Fc \deqs \{f_k :k\in \N\}$ of $H^1$-orthonormal functions satisfying
\beq
[f_k, h]_{1, 2} \eqs \tilde{\lambda}_k \, \ang{f_k, h}_2 \foral h\in \Hone
\eeq
with the values $\tilde{\lambda}_k$ being strictly positive. 

Define $\lamk : = \tilde{\lambda}_k - 1$ and let $\Lambda \deqs \{\lamk : k\in \N\}$.  
Then $f_k$ satisfies \eqref{e5.1} with Neumann Laplacian eigenvalue $\lambda = \lamk$, and consequently 
\beq \label{e5.5}
 \ang{f_j,f_k}_2  \eqs \delta_{jk} / (1 + \lamk) \xfor j, k\in \N. \eeq
Take $\Fctd \deqs \{(1 + \lambda_k)^{1/2}f_k : k\in \N\}$. 
These eigendata satisfy the following result.

\btm \label{T5.1}
Assume $(B1)$, $(B2)$ hold and $\Lambda, \ \Fc, \, \Fctd$ are defined as above.
Then $\Lambda$ is an increasing sequence with $\lamk \to \infty$, $\Fc$ is a maximal orthonormal set in $\Hone$ and $\Fctd$ is a maximal
orthonormal set in $\Lt$. \etm
\bpf
The bilinear form $a$ in \eqref{e5.3} automatically obeys condition (A1) of \cite{Au4}.
The form $m$ obeys conditions (A2) and (A4) of \cite{Au4} from our assumption (B2). 
Thus the first parts of this theorem follow since the $\tilde{\lambda}_k$ obey the results from Theorem 4.3 of \cite{Au4} . 

From Rellich's theorem  the imbedding of $\Hone$ into $\Lt$ is compact and also it has dense range, so 
condition (A5) holds and Theorem 4.6 of \cite{Au4} yields that $\Fctd$ is a maximal orthonormal set in $\Lt$.   \epf

Let $( \Lamo, \Fco)$ and $( \Lamt, \Fct)$ be sequences of Neumann Laplacian eigenvalues and eigenfunctions on $\Omo, \Omt$
generated as above with $\Lamo := \{ \lmok : k \in \N\}, \Fco := \{ \fok : k \in \N\} $ and similarly for $\Lamt, \Fct$.
Consider the family $\Uc$ consisting of the dyads
\beq
\ujk := \foj\otimes \ftk \xwith j, k\in \N.
\eeq
These dyads have the following properties.

\btm \label{T5.2}
Assume $\Omo, \Omt$ are regions that satisfy (B1) and (B2), $\Omp := \Omot$  and $\Lamo, \Lamt, \Fco, \Fct, \Uc$ as above. 
Then each $\ujk\in \Uc$ is a  Neumann Laplacian eigenfunction of \eqref{e5.1} on  $\Omp$ corresponding to the eigenvalue $\lmoj + \lmtk$.
 Moreover, $\Uc $ is orthogonal in $\Ltp$ and also in $\Honep$.
\etm 
\bpf  When  $h \in \Honep$, then $h(x, \cdot) \in \Hont$ for almost all $x \in \Omo$ and $h(\cdot,y) \in \Hono$ for almost all $y \in \Omt$.
From the Neumann eigenequation on $\Omo$,
\[ \IOmo \gradx \ujk(x,y) \cdot \gradx h(x,y) \ dx \eqs \lmoj \ftk(y) \IOmo  \foj(x) h(x,y) \ dx \]
for almost all $y \in \Omt.$
Integrating this over $\Omt$,  yields
\[   \IOmt \IOmo \, \ftk(y) \gradx \foj(x) \cdot \gradx h(x,y) \ dx \, dy \eqs \lmoj \IOmt \, \ftk(y) \IOmo \foj(x) h(x,y) \ dx \, dy. \]
Similarly one has that 
\[   \IOmo \IOmt \, \foj(x) \grady \ftk(y) \cdot \grady h(x,y) \ dy \, dx \eqs \lmtk \IOmo \,  \foj(x) \IOmt \ftk(y) h(x, y) \ dy \, dx. \]
Adding these shows that $\ujk$ is a solution of \eqref{e5.1} with eigenvalue $\lmoj + \lmtk$ as claimed.

The fact that the $\ujk$ are $L^2$-orthogonal follows by using Fubini's theorem and the fact that the families $\Fco, \Fct$ are 
$L^2$-orthogonal on $\Lto, \Ltt$ respectively.
This and the eigenequation \eqref{e5.1} yields the $H^1$-orthogonality of the $\ujk$.
\epf

From definition, the $\ujk$ satisfy
\beq \label{e5.9}
\ang{\ujk, \ujk}_{2, \Omp} \eqs \IOmo \foj(x)^2 \ dx \ \IOmt \ftk(y)^2 \ dy \eqs \frac{1}{(1 + \lmoj) \,(1 + \lmtk)}. \eeq
Thus, adding $\ang{\ujk, \ujk}_{2, \Omp}$ to both sides of \eqref{e5.1} gives
\beq \label{e5.7}
[\ujk, \ujk]_{1, 2,\Omp} \eqs  \frac{1 + \lmoj + \lmtk}{(1 + \lmoj)\ (1 + \lmtk)}.
\eeq  

Define $\vjk, \, \vjktd$ to be the functions
\beq
\vjk \deqs (1+\lmoj)^{1/2}(1 + \lmtk)^{1/2}\, \ujk \xand \vjktd \deqs \left( \frac{(1+\lmoj)(1 +\lmtk)}{1 + \lmoj + \lmtk} \right)^{1/2} \, \ujk
\eeq
and let $\Vc \deqs \{\vjk : j, k\in \N\}$ and $\Vctd \deqs \{\vjktd : j, k \in \N\}$.  
Then $\Vc$ is an orthonormal set in $\Honep$, $\Vctd$ is an orthonormal set in $\Ltp$ and the following holds.

\begin{cor} \label{C5.3}
Assume $\Omo, \Omt$ are regions that satisfy (B1) and (B2) and  $\Omp := \Omot$.
Then $\Vc$ is an orthonormal basis of $\Honep$, $\Vctd$ is an orthonormal basis of $\Ltp$ , and 
$\Honep = \Hono \otimes \Hont$.
\end{cor}
\bpf
Suppose $h \in \Honep$,  then from the eigenequation \eqref{e5.1} and the preceding theorem, 
\[ [h,\ujk]_{1, 2, \Omp} \eqs (1 + \lmoj + \lmtk) \ang{h, \ujk}_{2,\Omp} \foral j,k\in \N. \]
From \eqref{e5.7} each value $(1 + \lmoj + \lmtk)$ is strictly positive. 
If $\Vc$ is not a maximal orthogonal set in $\Honep$, then there is a nonzero $\hat{h} \in \Honep$ with $[\hat{h},\ujk]_{1, 2, \Omp} = 0$
for all $j,k$.
From the first line, $\ang{\hat{h}, \ujk}_{2, \Omp} = 0$. 
However $\Vctd$ is a maximal orthogonal set in $\Ltp$ since $\Fco, \Fct$ are maximal orthogonal sets in $\Lto, \Ltp$
respectively and $\Ltp = \Lto \otimes \Ltt$. 
Thus $\hat{h} \equiv 0$ in $\Ltp$ and a fortiori in $\Honep$.
Hence $\Vc$ is a maximal $H^1$-orthonormal set as claimed. 
Since $\Vc$ is an $H^1$-orthonormal basis of the Hilbert space $\Honep$  consisting of dyads of the form $\foj \otimes \ftk$, then 
$ \Honep = \Hono \otimes \Hont$ from Theorem \ref{TPt1}.
\epf

\vspace{1em}
%%%-----------------------------------------------------------------------------------------------------------------------------------------------------------------------------
%%%		SECTION:		Robin  Laplacian Eigenproblems on $\Om$ and $\Omp$.
%%%----------------------------------------------------------------------------------------------------------------------------------------------------------------------------

\section{Robin Laplacian Eigenproblems on $\Om$ and $\Omp$.} \label{s6}

The {\it  Robin Laplacian eigenproblem} on a region $\Om \subset \RN$ is the problem of finding nontrivial solutions
$(\lambda,u) \in \R \times \Hone$ of the system
\beq \label{e6.1}
\IOm \gradu \cdot \nabla h \, \, dx \pls \Iby b \, u \, h \dsg \eqs \lambda \IOm   \,  u \, h \, dx  \foral h\in \Hone.
\eeq
This is the weak form of the eigenvalue problem of finding nontrivial solutions of the system
\beq
-\Delta u(x)  \eqs  \lambda \,    u(x) \xon \Om  \xwith  \Dnu u + bu \eqs  0  \xon \bdy.
\eeq

A general theorem about the spectrum of second order Robin elliptic eigenproblems on bounded regions is 
given as Theorem 7.2 of \cite{Au4} and is based on a similar analysis to that of Theorem \ref{T5.1} above.
Let $\Lambda \deq \{\lamk : \, k \in \N \}$ be the sequence of eigenvalues of \eqref{e6.1} repeated according to multiplicity
and   $\Gc \deqs \{g_k : \, k \in \N\}$ be an associated sequence of Robin eigenfunctions normalized so that 
\beq \label{e6.2}
[g_j,g_k]_b \deqs \IOm \nabla g_j \cdot \nabla g_k \, dx \pls \Iby b\, g_j\, g_k \, \dsg \eqs \delta_{jk} \foral j,k \in \N. 
\eeq
Then the eigenequation \eqref{e6.1} implies  that 
\beq \label{e6.3}
  \ang{g_j,g_k}_2  \eqs \delta_{jk} / \lamk \xfor j,  k \geq  1 . \eeq

Define $\Gctd \deqs \{\lambda_k^{1/2}g_k : k\in \N\}$.  Then Theorem 7.2 of \cite{Au4}  may then be stated as
\btm \label{T6.1}
Assume that (B1)-(B3) hold and $\Lambda, \Gc, \Gctd $ are  defined as above. 
Then $\Lambda$ is an increasing sequence with $\lamk \to \infty$,  $\Gc$ is a maximal $b$-orthonormal set in $\Hone$ and $\Gctd$ is a maximal
orthonormal set in $\Lt$. \etm

To prove results about the Robin Laplacian eigenproblem on the product region $\Omp := \Omot$, first note
  that its boundary $\bdy_p$  is given by 
\[ \bdy_p \eqs ( \pal \Omo \times \Omb_2) \ \cup \ (\Omb_1 \times \ \pal \Omt). \]
Also, $\nu(x,y) \eqs (\nu_1(x), 0)$ for $(x,y) \in  \pal \Omo \times \Omt$ and similarly  $\nu(x,y) \eqs (0, \nu_2(y))$ for 
$(x,y) \in   \Omo \times \pal \Omt$. 
Here $\nu_1, \nu_2$ are unit normal vectors in $\Rno, \Rnt$ respectively.
The region $\Omp$ in general has ``corners" with no unit outward normal at points $(x,y) \in \bdy_1 \times  \bdy_2$.

The previous condition (B3) must be modified for Robin problems on $\Omp$. We now require that

\noindent   {\bf (B4):} \quad  {\it $b: \bdyp \to [0,\infty] $ is such that $b(x,y) = b_1(x)$ for $a.e. \; y \in \Omt$, and 
$b(x,y) = b_2(y)$ for $a.e. \; x \in \Omo$ with $b_1, b_2$ satisfying (B3) on $\bdyo, \bdyt$ respectively.}

Let $( \Lamo, \Gco)$ and $( \Lamt, \Gct)$ be the Robin Laplacian eigendata on $\Omo, \Omt$
generated as above with $\Lamo := \{ \lmoj : j \in \N\}, \Gc_1 := \{ \goj  : j\in \N \} $ and similarly for $\Lamt, \Gc_2$.
Consider the family $\Uc$ consisting of the dyads 
\beq
\ujk := \goj \otimes \gtk \xwith j, k\in \N.
\eeq
These functions have the following properties.

\btm \label{T6.2}
Assume that $(B1)$-$(B4)$ hold, $\Omp := \Omot$  and $\Lamo, \Lamt, \Gco, \Gct, \Uc $ as above. 
Then each $\ujk$ is a  Robin Laplacian eigenfunction of \eqref{e6.1} on  $\Omp$ corresponding to the eigenvalue $\lmoj + \lmtk$.
 Moreover $\Uc $ is an orthogonal family in $\Honep \ \wrt$ the $b$-inner product and is an orthogonal family 
 in $\Ltp$.
\etm 
\bpf
When $\goj, \gtk$ are as above then they are solutions, respectively, of
\begin{eqnarray} 
\IOmo \gradx \goj \cdot \gradx h_1 \, \, dx  \pls \Ibyo b_1 \, \goj \, h_1 \dsg_1 & \eqs & \lmoj  \ \IOmo   \,  \goj \, h_1 \, dx  \\
\IOmt \grady \gtk \cdot \grady  h_2 \, \, dy \pls \Ibyt b_2 \, \gtk \, h_2 \dsg _2 & \eqs &  \lmtk \ \IOmt   \,  \gtk \, h_2 \,  dy  
\end{eqnarray}
for all $h_1\in \Hono$ and $h_2\in \Hont$.
Multiply the first equation by $\gtk(y)h_2(y)$ and integrate over $\Omt$, and the second equation by $\goj(x)h_1(x)$ and integrate over $\Omo$.
Add these two equations and observe that the resulting  left hand side is  $[\ujk, h]_{b, \Omp}$ where $h(x,y) = h_1(x) \, h_2(y)$ as (B4) holds.
So this is $ [\ujk, h]_{b, \Omp} \eqs (\lmoj + \lmtk) \ \ang{\ujk, h}_{2, \Omp}$ for all such $h$. 
The class of such $h$ is dense in $\Honep$ from Corollary \ref{C5.3} above, so $\ujk$ is a Robin Laplacian eigenfunction
of \eqref{e6.1} on  $\Omp$ with eigenvalue $\lmoj + \lmtk$. 

Let $h_1 = \goj, h_2 = \gtk$ in the above formulae to obtain $[\goj, \goj]_{b, \Omo} = [\gtk, \gtk]_{b,\Omt} = 1$ and
\[ 1 \eqs \lmoj \ \ang{\goj, \goj}_{2, \Omo} \xand 1 = \lmtk\ \ang{\gtk, \gtk}_{2, \Omt} . \]
Multiply the first of these equalities by $\gtk^2$ and integrate over $\Omt$ and the second by $\goj^2$ and integrate over $\Omo$.
Then 
\beq \label{e6.7}
\frac{1}{\lmtk} \pls \frac{1}{\lmoj} \eqs (\lmoj + \lmtk) \ \IOmp \ujk(x,y)^2 \, dx \, dy \eqs [\ujk,\ujk]_{b,\Omp} \eeq

Choosing $h_1 = \goj$ with $j \neq k$ or $h_2 = \gtk$ with $k \neq j$ in these formulae shows that the functions $\ujk$ are orthogonal
$\wrt$ the $b$-inner product or the $L^2-$inner product on $\Omp$. \epf

Define $\wjk, \, \wjktd$ to be the functions
\beq
\wjk \deqs \left( \frac{\lmoj\, \lmtk}{\lmoj + \lmtk}\right)^{1/2} \ujk \xand \wjktd \deqs (\lmoj \, \lmtk)^{1/2} \ujk 
\eeq
and let $\Wc\deqs \{\wjk : j, k\in \N\}$ and $\Wctd \deqs \{ \wjktd : j, k \in \N\}$.  
Then $\Wc$ is a $b$-orthonormal set in $\Honep$ and $\Wctd$ is orthonormal in $\Ltp$.

\begin{cor} \label{C6.3}
Assume $(B1)$-$(B4)$ hold, $\Omp := \Omot$ and $\Wc, \, \Wctd$ as above.
Then $\Wc$ is a maximal $b$-orthonormall set in $\Honep$ and $ \Wctd$ a maximal orthonormal set in $\Ltp$  \end{cor}
\vspace{-0.5em}
\bpf
When $h \in \Honep$, then eigenequation \eqref{e6.1} and the preceding theorem give 
\[ [h, \ujk]_{b, \Omp} \eqs (\lmoj + \lmtk) \ang{h, \ujk}_{2,\Omp} \foral j,k\in \N. \]
If $\Wc$ is not a maximal $b$-orthonormal set in $\Honep$, 
then there is a nonzero $\hat{h} \in \Honep$ with $[\hat{h},\ujk]_{b,\Omp} = 0$ for all $j,k\in \N$.
Thus $\ang{\hat{h}, \ujk}_{2, \Omp} = 0$ for all $j,k$ as the individual Robin eigenvalues always are strictly positive.
However $\Gco, \Gct$ are maximal orthogonal sets in $\Lto, \Ltt$ respectively so their tensor  products are a maximal orthogonal set 
in $\Ltp$. Hence $\hat{h}=0$ in $\Ltp$ and thus also in $\Honep$. 
Therefore $\Wc$ is a maximal orthogonal set in $\Honep$. \epf

\vspace{2em}
%%%-----------------------------------------------------------------------------------------------------------------------------------------------------------------------------
%%%		SECTION:		D-N  Laplacian Eigenproblems on $ \Omp$.
%%%----------------------------------------------------------------------------------------------------------------------------------------------------------------------------

\section{Dirichlet-Neumann Laplacian Eigenproblems on $ \Omp$.} \label{s7}

In the preceding sections, the eigenproblems on the regions $\Omo, \Omt$ were of the same type.
In this section the Laplacian eigenproblem on $\Omp$ with Dirichlet conditions on $\Omo$ and Neumann
conditions on $\Omt$ is considered.

This problem can be posed in a number of different ways.
The simplest appears to be to seek eigenfunctions in the tensor product space $V := \Hzoo \otimes \Hont$.
% Need to define tensor products H_1 \otimes H_2 in section 2 - including tensor inner products 

It is straightforward to verify that \beq \label{e7.1} 
[u,v]_{\nabla, \Omp} \deqs \IOmp [\, \gradx u \cdot \gradx v \pls \grady u \cdot \grady v \, ] \, dx \, dy \eeq
is an  inner product to  this space V. 

The {\it  Dirichlet Neumann  Laplacian (DNL) eigenproblem} on the product region $\Omp$ is the problem of finding 
nontrivial solutions  $(\lambda,u) \in \R \times V$ of the system
\beq \label{e7.3}
\IOmp \gradu \cdot \nabla h \, \, dx \, dy  \eqs \lambda \IOmp    \,  u \, h \, dx \, dy  \foral h \in V.
\eeq
This is the weak form of the eigenvalue problem of finding nontrivial solutions of the usual eigenvalue equation
$\Delta u = \lambda u$ on $ \Omp$ subject to the  boundary conditions
\beq
  u \eqs 0 \xon \bdypo \xand  \Dnu u \eqs  0  \xon \bdypt
\eeq

Suppose now that $\Eco : = \{ \eoj : j\in \N\}$ is the family of Dirichlet Laplacian eigenfunctions on $\Omo$ defined as in Section \ref{s3}
and $\Fct := \{ \ftk : k \in \N \}$ is the family of Neumann Laplacian eigenfunctions on $\Omt$ defined as in Section \ref{s5}.
Let $\Lamo, \Lamt$ be the associated sequences of eigenvalues and consider the family  $\Uc$ of dyads 
\beq \label{e7.5}
\ujk \deqs \eoj \otimes \ftk \qquad \text{defined on } \Omot; \; j, k \in \N. \eeq
Each dyad $\ujk$ is in $V$ and the following holds.

\btm \label{T7.1}
Assume $(B1)$, $(B2)$ hold, and $\Lamo, \Lamt, \Eco, \Fct, \Uc$ as above. 
Then each $u_{jk} \in \Uc$ is a DNL eigenfunction of \eqref{e7.3} on  $\Omp$ corresponding to the eigenvalue $\lmoj + \lmtk$.
 Moreover  $\Uc $ is an orthogonal family in $V \ \wrt$ the inner product $[\cdot,\cdot]_{\nabla, \Omp}$ and is an orthogonal family 
 in $\Ltp$.
\etm 
\bpf
When  $h \in V$, then $h(x,\cdot) \in \Hont$ for $a.e. \; x \in \Omo$ and $h(\cdot,y) \in \Hzoo$ for $a.e. \; y \in \Omt$.
From the Dirichlet eigenequation on $\Omo$,
\[  \IOmo \, \ftk(y) \gradx \eoj(x) \cdot \gradx h(x,y) \ dx  \eqs \lmoj \ftk(y) \IOmo \eoj(x) \, h(x,y)  \ dx \]
for almost all $y \in \Omt.$
Integrating this over $\Omt$ then yields
\[   \IOmt\IOmo \gradx \ujk(x,y) \cdot \gradx h(x,y) \ dx \, dy \eqs \lmoj \IOmt  \IOmo \ujk(x,y)\,  h(x,y)\ dx \, dy. \]
Similarly one has that 
\[   \IOmo \IOmt \grady \ujk(x,y) \cdot \grady h(x,y) \ dy \, dx \eqs \lmtk \IOmo \,  \IOmt \ujk(x,y) \, h(x, y) \ dy \, dx. \]
Adding these shows that $\ujk$ is a solution of \eqref{e7.3} with eigenvalue $\lmoj + \lmtk$ as claimed.

The fact that the $\ujk$ are $\nabla$-orthogonal follows by using Fubini's theorem and the fact that the families $\Eco, \Fct$ are 
orthonormal  $\wrt$ the $\nabla$-inner product on $\Hzoo, \Hont$ respectively. 
These computations shows that the $\ujk$ have 
\beq \label{e7.7}
[\ujk, \ujk]_{\nabla,\Omp} \eqs \IOmo \IOmt  \, [\, |\gradx \ujk|^2 \pls  |\grady \ujk|^2 \, ] \ dy \, dx \eqs \frac{\lmoj + \lmtk}{1 + \lmtk}.
\eeq 
They are $L^2$-orthogonal on $\Omp$ since the individual factors are $L^2$-orthogonal on $\Omo, \Omt$ respectively. \epf

Note that the functions in $\Uc$ will be a maximal orthogonal set in $V$ from the general properties of tensor products of Hilbert spaces
as $\Eco, \Fct$ are maximal orthogonal sets in $\Hzoo, \Hont$ respectively.
Therefore, the above result implies that the eigenfunctions of this
DNL problem span $V$.

\vspace{2em}
%%%-----------------------------------------------------------------------------------------------------------------------------------------------------------------------------
%%%		SECTION:		D-R  Laplacian Eigenproblems on $ \Omp$.
%%%----------------------------------------------------------------------------------------------------------------------------------------------------------------------------

\section{Dirichlet-Robin Laplacian Eigenproblems on $ \Omp$.} \label{s8}

A similar analysis holds when Dirichlet conditions are imposed on $\Omo$ and Robin conditions are required on $\Omt$.
In this case  take $V = \Hzoo \otimes \Hont$ with  the inner product
\beq \label{e8.1} 
[u,v]_{b2} \deqs \IOmp  \gradu \cdot \gradv\, \,  dx \, dy  \pls \IOmo\Ibyt  b_2(y) u \, v \,  \dsg_2   \ dx \eeq
This is  easily verified to be an  inner product   when $b_2$ obeys (B3) on $\bdyt$. 

The {\it  Dirichlet Robin  Laplacian (DRL) eigenproblem} on the product region $\Omp$ is the problem of finding 
nontrivial solutions  $(\lambda,u) \in \R \times V$ of the system
\beq \label{e8.3}
[u,h]_{b2} \eqs     \lambda \ang{u,  h}_{2,\Omp}  \foral v \in V.    \eeq
This is the weak form of the eigenvalue problem of finding non-trivial solutions of the usual eigenvalue equation
$\Delta u = \lambda u$ on $ \Omp$ subject to the  boundary conditions
\beq
  u \eqs 0 \xon \bdypo \xand  \Dnu u \pls b_2 u  \eqs  0  \xon \bdypt
\eeq

Suppose $\Eco : = \{\eoj : j\in \N\}$ is the family of Dirichlet Laplacian eigenfunctions on $\Omo$ defined as in Section \ref{s3}
and that $\Gct := \{ \gtk : k \in \N \}$ is the family of Robin Laplacian eigenfunctions on $\Omt$ defined as in Section \ref{s6}.
Let $\Lamo, \Lamt$ be the associated sequences of eigenvalues and consider the family $\Uc$ of dyads  
\beq \label{e8.5}
\ujk \deqs \eoj\otimes \gtk  \xwith j, k \in \N. \eeq
Each of these dyads is in $V$ and the  following holds.

\btm \label{T8.1}
Assume $\Omo, \Omt, b_2$  satisfy $(B1)$-$(B3)$,  and $\Lamo, \Lamt, \Eco, \Gct,\Uc$ as above. 
Then each $\ujk  \in \Uc$ is a DRL eigenfunction of \eqref{e8.3} on  $\Omp$ corresponding to the eigenvalue $\lmoj + \lmtk$.
 Moreover  $\Uc $ is an orthogonal family in $V \ \wrt$ the inner product $[\cdot,\cdot]_{b2}$ and is an orthogonal family 
 in $\Ltp$.
\etm 
\bpf      
The fact that each $\ujk \in \Uc$ is a DRL eigenfunction of \eqref{e8.3} on  $\Omp$ with eigenvalue $\lmoj + \lmtk$
follows from the evaluation of the relevant integrals just as in the first part of the proof of Theorem \ref{T7.1}. 
This leads to $[u_{jk},h]_{b2} \eqs (\lmoj+\lmtk) \ang{\ujk,h}_{2, \Omp}$ for all $h \in V$.
Take $h = u_{lm}$ and use Fubini's theorem to verify that these functions are orthogonal in both the $L^2$- and $b2$-inner products.  
Note that a formula similar to \eqref{e7.7} can be obtained with the  $b2$-norm in place of $[\ujk, \ujk]_{\nabla, \Omp}$. 
\epf

\vspace{1em}
%%%-----------------------------------------------------------------------------------------------------------------------------------------------------------------------------
%%%		SECTION:		 Steklov  Laplacian Eigenproblems on $\Om \Omp$.
%%%----------------------------------------------------------------------------------------------------------------------------------------------------------------------------

\section{Dirichlet-Steklov, Neumann-Steklov, and Robin-Steklov Laplacian Eigenproblems on $ \Om$ and $\Omp$.} \label{s9}

The above sections all treated well-known eigenproblems for the Laplacian. 
Somewhat surprisingly, a similar result also holds for 2-parameter Dirichlet-Steklov (DSL), Neumann-Steklov (NSL)
and Robin-Steklov (RSL) eigenproblems for the Laplacian.  
These are problems that involve an eigenparameter $\lambda$ for the Laplacian and also a Steklov eigenvalue 
$\delta$ in the boundary condition on $\Omt$. 
Such problems arise in fluid mechanics and elsewhere, see Auchmuty and Simpkins \cite{AuS}.

The {\it weighted harmonic Steklov eigenproblem}  on a region $\Om \subset \RN$ is the problem of finding non-trivial
$(\delta, s) \in \R \times \Hone$ that satisfy the equation
\beq \label{e9.1}
[s,h]_{\nabla} \deqs \IOm \grads \cdot \nabla h \ dx \eqs \delta \ \Iby b \, s \, h \ \dsg \foral h \in \Hone.  \eeq 
Here $b$ is a function on $\bdy$  that obeys (B3). 
$b \equiv 1/|\bdy|$  is the standard harmonic Steklov eigenproblem.

This is the weak form of the problem of finding a solution of Laplace's equation on $\Om$ subject to
the boundary condition $\Dnu s \eqs \delta \, b \, s$ on $\bdy$. 
While this system was first studied by Steklov in 1902, the results used  here are taken from  
Auchmuty  \cite{Au3} and \cite{Au2} as well as Section 8 of \cite{Au4}.
Let $\Harm$ be the orthogonal complement of $\Hzone, \wrt$ the inner product of  \eqref{e4.2}.
Then  $\Harm$ is the space of all $H^1-$harmonic functions on $\Om$ when (B3) holds.

Observe that $\delta_0 = 0$ is an eigenvalue of \eqref{e9.1} with the associated eigenfunction $s_0(x) \equiv 1$.
This is the least eigenvalue and all other  Steklov eigenvalues must be strictly positive. 
Consider the sequence of constrained variational problems $\Pc_k$ as described by  (4.1) in \cite{Au4} taking the bilinear forms $a, m$ in that paper to be, respectively,
\beq
[u, v]_{b,\Om} \deqs \IOm \gradu\cdot\gradv dx \pls \Iby b\, u\, v\,  \dsg \xand \ang{u, v}_{b, \bdy} \deqs \Iby b\, u \, v \, \dsg
\eeq
This generates a sequence of eigenfunctions $s_k$ in $\Hone$ satisfying \eqref{e9.1} that are normalized with respect to the inner product $[\cdot, \cdot]_{b, \Om}$. 
Let $\Sc \deqs \{s_k : \, k \in \N \}$ be such a sequence and $\Lambda \deqs \{\delk: \, k \in \N \}$ be the associated sequence of eigenvalues.  
Then the eigenequation \eqref{e9.1} implies  that 
\beq
\begin{aligned}   \label{e9.5}
[s_k, s_k]_{\nabla}  & \eqs \frac{\delk}{1+\delk} & \xand  \ang{s_k, s_k}_{b,\bdy} & \eqs  \frac{1}{1+ \delk}   \foral k\in \N  \\
[ s_j, s_k ]_{\nabla}  & \eqs 0 & \xand   \ang{s_j \,  s_k}_{b, \bdy} & \eqs 0 \, \qquad \xwhen j \neq k. 
\end{aligned}
\eeq
The orthonormal set $\Scz := \Sc \cup \{s_o\}$ in $\Hone$ is now a basis of the subspace $\Harm$ of all $H^1$-harmonic 
functions on $\Om$ as stated below.

\btm \label{T9.1}
Assume that $(B1)$-$(B3)$ hold and $\Lambda, \Sc$ are  sequences defined by the successive variational problems $\Pc_k$.
Then $\Lambda$ is an increasing sequence with $\delk \to \infty$ and $ \Scz$ is a maximal orthogonal set in $\Harm$. \etm

This is a combination of Theorem 8.2 of \cite{Au4} and some analysis from \cite{Au2}. 
The later paper shows that under further conditions on $b$, the traces of the Steklov eigenfunctions on $\bdy$ are a maximal orthogonal set in 
$\Ltbby$

The {\it  Dirichlet-Steklov Laplacian (DSL) eigenproblem} on the product  region $\Omp := \Omot$ is the problem of finding nontrivial 
solutions $(\lambda, \delta,u) \in \R^2 \times V$ of the system
\beq \label{e9.3}
\IOmp \gradu \cdot \nabla h \, \, dx \, dy   \eqs \lambda \IOmp     u \, h \, dx \, dy \pls \delta \, \IOmo \Ibyt  b_2 \, u \, h \dsg_2 \, dx  \foral h \in V.
\eeq
Here $V := \Hzoo \otimes \Hont$ as in Section \ref{s7}, $b_2$ is a function on $\bdyt$ that obeys (B3), and the left hand side of \eqref{e9.3} is  $[u,v]_{\nabla, \Omp}$.

 This is  a 2-parameter eigenproblem and is the weak form of the problem of finding nontrivial solutions of the system 
 $-\Delta u  =  \lambda \, u$ on $\Omp$   subject to the boundary conditions
\beq \label{e9.4}
   u  \eqs  0  \xon \bdyo \times \Omt \xand \Dnu u \eqs \delta \, b_2 \, u \xon \Omo \times \bdyt. 
\eeq

Suppose now that $\Eco :=\{\eoj : j\in \N\}$ and $\Lamo := \{\lmoj : j\in \N\}$ is the Dirichlet Laplacian eigendata on $\Omo$ defined as in Section 3,
and let $\Sct := \{ \stk : k \in\N_0 \}$ and $\Lamt : = \{\deltk : k\in \N_0\}$ be the harmonic Steklov eigendata on $\Omt$ defined as above, 
with $\Omt$ in place of $\Om$ and $\N_0 : = \N\cup\{0\}$.
Let $\Uc$ be the family consisting of the dyads  
\beq \label{e9.7}
\ujk \deqs \eoj \otimes \stk  \xwith j \in \N, \, k \in \N_0. \eeq
Then each of these dyads is in $V$ and the following holds.

\btm \label{T9.2}
Assume $\Omo, \Omt$ are regions that satisfy $(B1)$ and $(B2)$, $b_2$ satisfies $(B3)$ on $\bdyt$, and define $\Lamo, \Lamt, \Eco, \Sct, \Uc $ as above.
Then each $u_{jk} \in \Uc$ is a DSL eigenfunction of \eqref{e9.3} on  $\Omp$ corresponding to the eigenpair $(\lmoj, \deltk  )$. \etm 
\bpf
This result follows from a straightforward computation and the use of Fubini's theorem.
When  $h \in V$, then $h(x,\cdot) \in \Hont$ for $a.e.\; x \in \Omo$ and $h(\cdot,y) \in \Hzoo$ for $a.e.\; y \in \Omt$.
From the Dirichlet eigenequation on $\Omo$,
\[ \IOmo \stk(y) \gradx \eoj(x) \cdot \gradx h(x,y) \ dx \eqs \lmoj \stk(y) \IOmo \eoj(x)\, h(x,y) \ dx \]
for almost all $y \in \Omt.$
Integrating this over $\Omt$, then yields
\[   \IOmp  \gradx \ujk(x,y) \cdot \gradx h(x,y) \ dx \, dy \eqs \lmoj \IOmp \ujk(x,y) \, h(x,y)\ dx \, dy. \]
Similarly, from the Steklov eigenequation on $\Omt$ one has that 
\[   \IOmp  \grady \ujk(x,y) \cdot \grady h(x,y) \ dy \, dx \eqs \deltk \IOmo \,  \Ibyt b_2(y) \,  \ujk(x,y) \, h(x,y) \ \dsg_2 \, dx. \]
Adding these shows that $\ujk$ is a solution of \eqref{e9.3} with eigenpair $(\lmoj, \deltk)$ as claimed.
\epf

From the preceding Theorems \ref{T3.1} and \ref{T9.1}, the sets $\Eco$ and $\Sct$ are maximal orthonormal sets in $\Hzone$ and $\Harm$
so they will be a maximal linearly independent set in the space $\Hzone \otimes \Harm$. 
Unfortunately the  computation of $[u_{jk}, u_{lm}]_{\nabla, \Omp}$ with the formulae from the above proof, shows that this inner product on $V$
may be nonzero when $j=l$ and $k \neq m$ since one has that $s_{2k}, s_{2m}$ are not, in general, $L^2$-orthogonal on $\Omt$.

The {\it  Neumann-Steklov Laplacian (NSL) eigenproblem} on the product  region $\Omp := \Omot$ is the problem of finding nontrivial 
solutions $(\lambda, \delta,u) \in \R^2 \times \Honep$ of the system
\beq \label{e9.8}
\IOmp \gradu \cdot \nabla h \, \, dx \, dy   \eqs \lambda \IOmp     u \, h \, dx \, dy \pls \delta \, \IOmo \Ibyt  b_2 \, u \, h \dsg_2 \, dx.
\eeq
for all $h\in \Hone$.
This is the weak form of the eigenvalue problem of finding nontrivial solutions of the system $-\Delta u =  \lambda \, u $ on $\Omp$
 subject to the boundary conditions
\beq \label{e9.9}
  \Dnu u  \eqs  0  \xon \bdyo \times \Omt \xand \Dnu u \eqs \delta \, b_2 \, u \xon \Omo \times \bdyt. 
\eeq

Suppose now that $\Fco:=\{\foj : j\in \N\}$ and $\Lamo :\{\lmoj : j\in \N\}$ is the Neumann Laplacian eigendata on $\Omo$ defined as in Section \ref{s5},
and $\Sct, \Lamt$ is the harmonic Steklov eigendata on $\Omt$ defined as above.
Let $\Uc$ be the family consisting of the dyads 
\beq \label{e9.11}
\ujk  \deqs \foj  \otimes \stk  \xwith j\in \N,\;  k \in \N_0, \eeq
Each of these dyads is in $\Honep$ and the  following holds.

\btm \label{T9.3}
Assume $\Omo, \Omt$ are regions that satisfy $(B1)$ and $(B2)$, $b_2$ satisfies $(B3)$ on $\bdyt$, and define $\Lamo, \Lamt, \Fco, \Sct, \Uc $ as above.
Then each $u_{jk} \in \Uc$ is an NSL eigenfunction of \eqref{e9.8} on  $\Omp$ corresponding to the eigenpair $(\lmoj,  \delk)$. \etm 
\bpf
This result follows from a straightforward computation and the use of Fubini's theorem.
When  $h \in \Honep$, then $h(x,\cdot) \in \Hont$ for $a.e. \; x \in \Omo$ and $h(\cdot,y) \in \Hono$ for $a.e. \; y \in \Omt$.
From the Neumann eigenequation on $\Omo$,
\[ \IOmo \stk(y) \gradx \foj(x) \cdot \gradx h(x,y) \ dx \eqs \lmoj \stk(y) \IOmo \foj(x) \, h(x,y) \ dx \]
for almost all $y \in \Omt.$
Integrating this over $\Omt$, then yields
\[   \IOmp  \gradx \ujk(x,y) \cdot \gradx h(x,y) \ dx \, dy \eqs \lmoj \IOmp \ujk(x,y) \, v(x,y)\ dx \, dy. \]
Similarly one has that 
\[   \IOmp \grady \ujk(x,y) \cdot \grady h(x,y) \ dy \, dx \eqs \deltk \IOmo \,  \Ibyt b_2(y) \,  \ujk(x,y) \, h(x,y) \ \dsg_2 \, dx. \]
Adding these shows that $\ujk$ is a solution of \eqref{e9.3} with eigenpair $(\lmoj, \deltk)$ as claimed.
\epf

From the preceding Theorems \ref{T5.1} and \ref{T9.1}, the sets $\Fco$ and $\Sct$ are maximal orthonormal sets in $\Hono$ and $\Hc(\Omt)$
so they will be a maximal linearly independent set in the space $\Hone \otimes \Harm$. 
Again inner products of the form $[u_{jk}, u_{lm}]_{\nabla, \Omp}$ need not always be zero when $k \neq m$ so there are only partial
 orthogonality results in this case.

The {\it Robin-Steklov Laplacian (RSL) eigenproblem} on the product region $\Omp \deqs \Omo\times \Omt$ is the problem of finding nontrivial solutions $(\lambda, \delta, u) \in \R^2 \times \Honep$ of the system
\beq\label{e.robin}
\IOmp \gradu \cdot \nabla h \, dx\, dy + \IOmt\Ibyo b_1 \ujk h\dsg_1 dy \, = \,  \lambda \IOmp u\, h\, dx\, dy + \delta \IOmo\Ibyt b_2 u\, h \, \dsg_2 dx
\eeq

This is the weak form of the eigenvalue problem of finding nontrivial solutions of the system $-\Delta u = \lambda u$ on $\Omp$ subject to the boundary conditions
\beq
\Dnu u + b_1 u = 0 \xon \bdyo\times \Omt \xand \Dnu u = \delta b_2 u \xon \Omo\times \bdyt.
\eeq

Suppose now that $\Gco := \{\goj : j\in \N\}$ and $\Lamo : =\{\lmoj : j\in \N\}$ is the Robin Laplacian eigendata on $\Omo$ defined as in Section \ref{s6}, 
and $\Sct, \Lamt$ is the harmonic Steklov eigendata defined as above.  
Let $\Uc$ be the family consisting of the dyads
\beq
\ujk \deqs \goj\otimes \stk \xwith j\in \N, k\in \N_0.
\eeq
Each of these dyads is in $\Honep$ and the following holds.

\btm\label{T9.4}
Assume $\Omo, \Omt$ are regions that satisfy $(B1)$ and $(B2)$, $b$ satisfies $(B4)$ on $\partial \Omp$, 
and $\Lamo, \Lamt, \Gco, \Sct, \Uc$ as above.
Then each $\ujk \in \Uc$ is a $RSL$ eigenfunction of \eqref{e.robin} corresponding to the eigenpair $(\lmoj, \deltk)$.
\etm
\bpf
As in the previous theorems, this result follows from a straightforward computation and the use of Fubini's theorem.
When $h\in \Honep$, then $h(x, \cdot)\in \Hont$ for $a.e. \; x\in \Omo$ and $h(\cdot, y)\in \Hono$ for $a.e. \; y\in \Omt$.
The Robin eigenequation on $\Omo$ holds for $u= \goj$, $\lambda = \lmoj$ and test function $h(\cdot,y)$ for almost every $y\in \Omt$.
Multiply that equation by $\stk$ and integrate over $\Omt$ to obtain
\beq
\IOmp \gradx \ujk \cdot \gradx h \, dx\, dy + \IOmt\Ibyo b_1 \ujk h\dsg_1 dy \, = \,  \lmoj \IOmp \ujk\, h\, dx\, dy .
\eeq
Similarly from the Steklov eigenequation on $\Omt$ one obtains
\beq
\IOmp \grady \ujk \cdot \grady h \, dx\, dy   \, = \,  \deltk \IOmo\Ibyt b_2 \ujk \, h \, \dsg_2 dx.
\eeq
Adding these shows that $\ujk$ is a solution of \eqref{e.robin} with eigenpair $(\lmoj, \deltk)$ as claimed.
\epf

From the preceding Theorems \ref{T6.1} and \ref{T9.1}, the sets $\Gco$, and $\Sct$ are maximal orthonormal sets in $\Hono$ and $\Hc(\Omt)$ so they will be a maximal linearly independent set in the space $\Hono\otimes \Hc(\Omt)$.

% \vspace{1em}
%%%-----------------------------------------------------------------------------------------------------------------------------------------------------------------------------
%%%		SECTION:		Boundary Trace Spaces of $ \Omp$.
%%%----------------------------------------------------------------------------------------------------------------------------------------------------------------------------

\section{Boundary Trace Spaces for  $\Omp$.} \label{s11}

A well-known characterization of the trace space $\Hhby$ of a region $\Om \subset \RN$ is that it is isomorphic
to the quotient space  $\Hone / \Hzone$ or, equivalently, to the orthogonal complement of $\Hzone$ in $\Hone$. 
Here we shall describe this orthogonal complement when $\Om$ is a product region and, in particular, obtain  orthogonal bases of  this space.

The definition of the space $\Harm$ implies that 
\beq \label{e11.1}
\Hone \eqs \Hzone \oplusb \Harm \eeq 
when $\Om$ satisfies (B1) and (B3) holds. 
Thus the harmonic Steklov eigenfunctions on a region provide a basis for the trace space of $H^1-$functions. 
See Auchmuty \cite{Au2} for an analysis of this situation.
In view of this one can ask about possible decompositions of the space of harmonic functions on $\Omp$ in terms of
functions on $\Omo, \Omt$.

A simple calculation shows that the tensor product of two Steklov eigenfunctions on regions need not be a Steklov eigenfunction on  the product region.
Also  there are Steklov eigenfunctions on product regions that are not obtained as finite linear combinations
of  Steklov eigenfunctions on the individual components.
When $\Om \subset \R$, then there are only two Steklov eigenfunctions of the Laplacian on $\Om$ while there are infinitely many Steklov  eigenfunctions on a box. 
See Auchmuty and Cho \cite{AuC}, or Girouard and Polterovich \cite{GP} for explicit  formulae  for the harmonic Steklov eigenfunctions  on boxes in  the plane.

The preceding results lead to the following isomorphism for  the space of harmonic functions on $\Omp$ when
$b \equiv 1$ on $\pal \Om$.

\btm \label{T11.1}
Suppose that $\Omp := \Omo \times \Omt$ and both $\Omo, \Omt$ satisfy $(B2), \ b \equiv 1$ on $\pal \Omp$
and $\Honep$ has this b-inner product.
Then $\Hc(\Omp)$ is linearly isomorphic to $V_0 \oplus V_1 \oplus V_2$ where 
\beq  \label{e11.3}
V_0 \eqs \Hc(\Omo) \otimes  \Hc(\Omt), \ \  V_1 \eqs H_0^1(\Omo) \otimes  \Hc(\Omt), \ \ 
V_2 \eqs   \Hc(\Omo)  \otimes H_0^1(\Omt) . \eeq 
Moreover if $\Sc_1, \Sc_2$ are orthonormal bases of $\Hc(\Omo), \Hc(\Omt)$ respectively,
and $\Eco, \Ect$ are orthonormal bases of $H_0^1(\Omo), H_0^1(\Omt)$ as in Section \ref{s3}, then
the sets $\Scott, \; \Eco \otimes \Sct, \;  \Sco \otimes \Ect, $ are maximal orthogonal sets in $V_0, V_1, V_2$
respectively.  \etm     
\bpf
From Theorem \ref{TPt2} upon using \eqref{e11.1} for $\Omo$ and $\Omt$ as well as \eqref{e11.1}
for the region $\Omp$, one sees that
\[ \Honep \eqs \Hzonep \oplus V \xwith V = V_0 \oplus V_1 \oplus V_2. \]
with  this b-inner product being used.
Thus $V$ and $\Hc(\Omp)$ are linearly isomorphic as they are orthogonal complements of $\Hzonep$ in
$\Honep$.
The characterization of the  bases of $V_0, V_1, V_2$ follows from the properties of bases of Hilbert tensor products.
\epf

This construction of an orthonormal  basis for $\Hc(\Omp)$ may be used to provide a different description of the
boundary traces of $H^1$-functions on $\Omp$.
It implies that  any function with nonzero boundary trace can be uniquely written as a function in $V$ and 
one with zero trace.     
Thus the $H^1$-boundary trace space is linearly isomorphic to $V$.

%   line 1230 approx
% ########    Bibliography  #########
%


\vspace{2em}


\begin{thebibliography}{AMR}

\vspace{1em}
\bibitem{Amick} C.J. Amick, ``Some remarks on Rellich's theorem and the Poincar\'{e} inequality", 
J. London Math. Soc. (2) {\bf 18 } (1973), 81 - 93.

\vspace{1em}
\bibitem{ABM} H. Attouch, G. Buttazzo and G. Michaille, 
\emph{Variational Analysis in Sobolev and BV Spaces}
SIAM, Philadelphia, (2006). 

\vspace{1em}
\bibitem{Au3} G. Auchmuty, ``Steklov eigenproblems and the representation of solutions of elliptic boundary value problems",
Numerical Functional Analysis and Optimization, {\bf 25 } (2004) 321-348.

\vspace{1em}
\bibitem{Au2} G. Auchmuty, ``Spectral Characterization of the Trace Spaces $\Hsby$",  
SIAM J of Mathematical Analysis, {\bf 38} (2006), 894-907.

\vspace{1em}
\bibitem{Au4}    G. Auchmuty, 
``Bases and Comparison Results for Linear Elliptic Eigenproblems",
J.Math. Analysis and Applications \textbf{390} (2012), pp. 394-406

\vspace{1em}   \bibitem{AuS}  G. Auchmuty and D.R. Simpkins, 
``Spectral representations and Approximations of Divergence-free Vector Fields", 
Q. of Applied Mathematics, (to appear). 

\vspace{1em}
\bibitem{AuC}  G. Auchmuty and M. Cho, 
``Boundary Integrals and Approximations of Harmonic Functions", 
Numer. Functional Analysis and Optimization, {.}, (2015).  

\vspace{1em}
\bibitem{BB1} P. Blanchard and E. Br\"{u}ning,
\emph{Variational Methods in Mathematical Physics},
Springer-Verlag, Berlin (1992).

\vspace{1em}
\bibitem{EG}  L.C. Evans and R.F. Gariepy,   
\emph{Measure Theory and Fine Properties of Functions},   CRC Press, Boca Raton (1992).

\vspace{1em}
\bibitem{GP}  A. Girouard and I. Polterovich, ``Spectral Geometry of the Steklov Problem", arXiv:1411.6567v1.

\vspace{1em}
\bibitem{St} W.A. Strauss,
\emph{Partial Differential Equations : An Introduction},
Second Edition, Wiley, New Jersey (2008). 

\vspace{1em}
\bibitem{Z2}  E. Zeidler,   \emph{Nonlinear Functional Analysis and its Applications, III, Variational Methods and Optimization}, 
Springer Verlag, New York (1985).



\end{thebibliography}
\end{document}